\documentclass{amsart}
\usepackage{amsmath}
\usepackage{amssymb}
\usepackage{graphicx,epsf,amsmath}
\usepackage{color}

\newtheorem{theorem}{Theorem}[section]

\newtheorem*{theorem*}{Theorem}

\theoremstyle{definition}

\newcommand{\inn}{{\quad\hbox{in } }}

\newcommand{\pp}{ {\partial} }

\newcommand{\R} {\mathbb R}
\newcommand{\cuad}{{\sqcap\kern-.68em\sqcup}}

\newcommand{\equ}[1]{(\ref{#1})}

\newcommand{\be}{\begin{equation}}
\newcommand{\ee}{\end{equation}}
\newcommand{\bee}{\begin{equation*}}
\newcommand{\eee}{\end{equation*}}
\newcommand{\bea}{\begin{eqnarray}}
\newcommand{\eea}{\end{eqnarray}}
\newcommand{\bs}{\begin{split}}
\newcommand{\es}{\end{split}}

\numberwithin{equation}{section}

\title[]{Non degeneracy of the bubble in the critical case for non local equations}

\author[J. Davila]{Juan D\'avila}
\address{\noindent J. D\'avila- Departamento de
Ingenier\'{\i}a  Matem\'atica and CMM, Universidad de Chile, Casilla
170 Correo 3, Santiago, Chile.} \email{jdavila@dim.uchile.cl}

\author[M. del Pino]{Manuel del Pino}
\address{\noindent M. del Pino- Departamento de
Ingenier\'{\i}a  Matem\'atica and CMM, Universidad de Chile, Casilla
170 Correo 3, Santiago, Chile.} \email{delpino@dim.uchile.cl}

\author[Y. Sire]{Yannick Sire}
\address{Y. Sire- Universit\'e Aix-Marseille 3, Paul C\'ezanne -- LATP  Marseille, France.}
\email{sire@cmi.univ-mrs.fr}

\begin{document}
\begin{abstract}
We prove the nondegeneracy of the extremals of the fractional Sobolev inequality as  solutions of a critical semilinear 
nonlocal equation  involving the fractional Laplacian.
\end{abstract}

\maketitle

\section{Introduction}

This paper establishes the {\em linear non-degeneracy} property of the extremals of the optimal Hardy-Littlewood-Sobolev inequality, which states the existence of a positive number $S$ such that for all $u\in C_0^\infty (\R^N)$ one has
\begin{equation}
 S\|u\|_{L^{2^*}(\R^N)}\ \le \   \|(-\Delta)^{\frac s2}u\|_{L^2(\R^N)} .
\label{sob}\end{equation}
Here $0< s< 1$, $N>2s$ and $2^* = \frac{2N}{N-2s}$.

The validity of this inequality, without its optimal constant, traces back to Hardy and Littlewood \cite{hl1,hl2}  and Sobolev \cite{sobolev}.  For $s=1$, \equ{sob} corresponds to the classical Sobolev inequality
 \begin{equation}
 S\|u\|_{L^{\frac{2N}{N-2}}(\R^N)}\ \le \   \|\nabla u\|^2_{L^2(\R^N)} .
\label{sob1}\end{equation}
Aubin \cite{aubin} and Talenti \cite{talenti} found the optimal constant and extremals for inequality \equ{sob1}. Indeed, equality
is achieved precisely by the functions
\be
w_{\mu,\xi}(x)  =   \alpha \left (   \frac{ \mu  } { \mu^2 + |x-\xi|^2}
\right )^{\frac{N-2}2}
\label{wmu}\ee
which for the  choice $\alpha= (N(N-2))^{\frac{N-2}4} $  solve
 the  equation
\be -\Delta w  = w^{\frac{N+2}{N-2}} , \quad w>0 \inn \R^N.  \label{w}\ee
The solutions \equ{wmu} are indeed the only ones of \equ{w}, see Caffarelli, Gidas and Spruk \cite{cgs}.

\medskip
For $s\ne 1$, the optimal constant in inequality \equ{sob} was first found by Lieb in 1983 \cite{lieb}, while alternative proofs have been provided   by by Carlen and Loss \cite{cl} and Frank and Lieb \cite{frank-lieb,fl2}.
Lieb established that the extremals correspond precisely to functions of the form
\be
w_{\mu,\xi}(x)  =   \alpha \left (   \frac{ \mu  } { \mu^2 + |x-\xi|^2}
\right )^{\frac{N-2s}2}, \quad \alpha >0,
\label{wmu2}\ee
which for a suitable choice  $\alpha = \alpha_{N,s}$ solve the equation
\be ( -\Delta)^s w  = w^{\frac{N+2s}{N-2s}}, \quad w>0 \inn \R^N.  \label{ws}\ee
Besides, under suitable decay assumptions, these are the only solutions of the equation, see
Chen, Li and  Ou \cite{chen-etal}, Li \cite{yanyan}, and  Li and Zhu \cite{lz}.

\medskip 
In equation \equ{ws} and in what remains of this paper, we will always mean that 
the equation
$$(-\Delta)^s u = f
\quad
\text{in } \R^N,
$$ 
is satisfied if 
\begin{align}
\label{inverse}
u(x)=  (-\Delta)^{-s} f (x)\, := \, \gamma\int_{\R^N}\frac{f(y)}{|x-y|^{N-2s}}\,dy,
\end{align}
as long as  $f$  has enough decay for the integral to be well defined.
The constant $\gamma=\gamma_{N,s}>0$ is so that \eqref{inverse} defines an inverse for the operator whose Fourier multiplier is  $|\xi|^{2s}$, namely,
$$
|\xi|^{2 s}\hat u(\xi) = \hat f(\xi).
$$
There are other definitions of $(-\Delta)^s u$,
which are equivalent to \eqref{inverse} under suitable assumptions,
 see \cite{caffarelli-silvestre}.


\medskip
When analyzing bubbling (blowing-up) behavior of semilinear elliptic equations involving critical Sobolev growth such as
Yamabe type problems, a crucial ingredient is understanding the {\em linear nondegeneracy} of the solutions \equ{wmu} of equation
\equ{w}. Let us observe that
\be
w_{\mu,\xi}(x)  =   \mu^{\frac{N-2}2 } \, w( \mu(x-\xi) ), \quad w(x) =  \alpha_N \left (   \frac{ 1  } { 1 + |x|^2}
\right )^{\frac{N-2}2}
\label{wmu1}\ee
which actually reflects the invariance of the equation under the above scaling and translations.
We say that  the solution $w$  is non-degenerate if  the generators of these operations are the only nontrivial elements in
the kernel of the associated linearized operator.  This property of Equation  \equ{w} is well-known, see for instance \cite{rey}.

\medskip
The purpose of this note is to establish the validity of this property for the solutions \equ{wmu2} of the fractional critical equation
\equ{ws}.  We observe now that
\be
w_{\mu,\xi}(x)  =   \mu^{\frac{N-2s}2 } \, w( \mu(x-\xi) ), \quad w(x) =  \alpha_{N,s} \left (   \frac{ 1  } { 1 + |x|^2}
\right )^{\frac{N-2s}2} .
\label{wmus}\ee
Differentiating  the equation 
\be ( -\Delta)^s w_{\mu,\xi}  = w_{\mu,\xi}^p  \inn \R^N, \quad p= {\frac{N+2s}{N-2s}}  \nonumber\ee
with respect of the parameters at $\mu=1$, $\xi=0$
we see that  the functions 
\be
\pp_{\mu} w_{\mu, \xi}   \ =\ \frac{N-2s}{2} w + x \cdot \nabla w, \quad
\pp_{\xi_i} w_{\mu, \xi}   \ =\ -\pp_{x_i} w 
\label{pp}\ee
annihilate the linearized operator around $w$, namely they satisfy the equation  
\begin{equation}\label{yamabe}
(-\Delta)^s \phi  =  p w^{p-1} \phi
\quad\text{in } \R^N
\end{equation}
 Our main result states the nondegeneracy of the solution $w$ in the following sense.

\begin{theorem}
\label{lnd}
The solution  $$w(x) =  \alpha_{N,s} \left (   \frac{ 1  } { 1 + |x|^2}
\right )^{\frac{N-2s}2}$$  of equation \equ{wmu2} is nondegenerate in the sense that all bounded solutions of equation \equ{yamabe} 
are linear combinations of the functions 
\begin{align}
\label{kernel}
\frac{N-2s}{2} w+ x \cdot \nabla w
,\quad\quad\quad
\pp_{x_i}w , \quad 1 \le i \le N  .
\end{align}

\end{theorem}

\section{Proof of Theorem~\ref{lnd}}

Let  $\phi$ be a bounded solution of equation \eqref{yamabe}.  Then $\phi$ solves the integral equation
\begin{align}
\label{formulation}
\phi(x)  =
\gamma \int_{\R^N} \frac{p w(y)^{p-1} \phi(y) }{|x-y|^{N-2s}} \, d y , \quad\forall x \in \R^N,
\end{align}
as we note that the integral is well defined thanks to the decay of $w$.

\medskip
The proof of the theorem is based on transforming equation \equ{formulation} 
into the sphere by means of stereographic projection. 
Let $\mathbb S^N$ be the sphere (contained in $\R^{N+1}$) and
$S : \R^N \to \mathbb S^{N} \setminus \{ pole \}$,
$$
S(x) =
\left(
\frac{2x}{1+|x|^2} , \frac{1-|x|^2}{1+|x|^2}
\right)
$$
be the stereographic projection. We have the following well-known formula for the Jacobian of the stereograhic projection:
$$
J_S(x) =  ( \frac{2}{1+|x|^2} ) ^N
= c w(x)^{\frac{N-2s}{2 N}} = c w(x)^{\frac{1}{p^* +1}} .
$$
Here and later $c>0$ is a constant that depends on $N,s$.
For $\varphi :\R^N \to \R$ define $\tilde\varphi :\mathbb S^N \to \R$ by
\begin{align}
\label{d1}
\varphi(x) = J(x)^{\frac{N+2s}{2 N}} \tilde \varphi( S(x)) .
\end{align}
Then we have
\begin{align}
\label{id1}
\int_{\R^N}
\int_{\R^N}
\frac{\varphi(x) \psi(y)}{|x-y|^{N-2s}}
\, d x \, d y
= \int_{\mathbb S^N}
\int_{\mathbb S^N} \frac{\tilde \varphi(\omega) \tilde \psi(\eta)}{|\omega-\eta|^{N-2s}} \, d \omega \, d \eta ,
\end{align}
where $|\omega-\eta|$ is the Euclidean distance in $\R^{N+1}$ from $\omega$ to $\eta$, see \cite{frank-lieb}.

\medskip

We rewrite  \eqref{formulation}
as
\begin{align}
\label{eq phi}
\phi(x) &
=
c \int_{\R^N} \frac{J(y)^{\frac{2 s}{N}} \phi(y) }{|x-y|^{N-2s}} \, d y  ,\quad\text{for all } x \in \R^N.
\end{align}
Let $\psi \in C_c^\infty(\R^N)$. Multiplying equation \eqref{eq phi} by $\psi J^{-\frac{2 s}{N}}$ and integrating in $\R^N$, we get
$$
\int_{\R^N} J(x)^{-\frac{2s}{N} } \phi(x) \psi(x) \, d x
=
\int_{\R^N}
\int_{\R^N}
\frac{J(y)^{\frac{2 s}{N}} \phi(y) J(x)^{-\frac{2 s}{N}} \psi(x)}{|x-y|^{N-2s}} \, d y \, d x .
$$
Let  $\tilde \phi$, $\tilde \psi$ be defined by \eqref{d1}. Then,  using \eqref{id1} we obtain
$$
\int_{\R^N} J(x)^{-\frac{2s}{N} } \phi(x) \psi(x) \, d x
=
\int_{\mathbb S^N}
\int_{\mathbb S^N}
\frac{J(S^{-1}(\omega))^{\frac{2 s}{N}} \tilde \phi(\omega)
J(S^{-1}(\eta) )^{-\frac{2 s}{N}} \tilde \psi(\eta)}{|\omega-\eta|^{N-2s}} \, d \omega \, d \eta.
$$
But
$$
\int_{\R^N} J(x)^{-\frac{2s}{N} } \phi(x) \psi(x) \, d x
=
\int_{\R^N} J(x) \tilde \phi(S(x) ) \tilde \psi(S(x)) \, d x
=
\int_{\mathbb S^N}
\tilde \phi(\omega) \tilde \psi(\omega) d \, \omega,
$$
so, replacing $\tilde \psi$ by a new test-function of the form
$\tilde \psi J^{\frac{2s}{N}}\circ S^{-1}$
we get
$$
\int_{\mathbb S^N}
J^{\frac{2s}{N}}\circ S^{-1}\tilde \phi \tilde \psi
=
c
\int_{\mathbb S^N}
\int_{\mathbb S^N}
\frac{J^{\frac{2s}{N}}\circ S^{-1}(\omega) \tilde \phi(\omega) \tilde \psi(\eta)}
{ |\omega-\eta|^{N-2s}}
\, d \omega
\, d \eta .
$$
Since $\tilde \psi$ is arbitrary (smooth with support away from the pole)
the function
\begin{align}
\label{def h}
h= J^{\frac{2s}{N}}\circ S^{-1}\tilde \phi
\end{align}
satisfies the integral equation
\begin{align}
\label{eh}
a h(\omega) =  \int_{\mathbb S^N} \frac{h(\eta)}{|\omega-\eta|^{N-2s}} \, d \eta  ,
\end{align}
for all $\omega \not=pole$, where $a>0$ depends only on $N,s$.
Note that  $h$ is a bounded function.
Indeed, in terms of $\phi$ this is equivalent to the estimate
\begin{align}
\label{est phi}
|\phi(x)|\le \frac{C}{|x|^{N-2s} } \quad \text{for all } |x|\ge 1.
\end{align}
We observe that if $v:\R^N \to \R$ satisfies $|v(x)| \le C/(1+|x|)^\nu$ and $\nu \ge 0$, then
\begin{align}
\label{conv}
\left|
\int_{\R^N} \frac{J(y)^{\frac{2 s}{N}}  }{|x-y|^{N-2s}} v(y) \, d y
\right|
\le \frac{C}{(1+|x|)^{\min(\nu +2 s,N-2s)}}.
\end{align}
Hence, starting from the hypothesis that $\phi $ is bounded, we get estimate
\eqref{est phi} after a finite number of applications of \eqref{conv}.
This implies that $h$ is bounded.
Then equation \eqref{eh} shows that $h$ has a continuous extension to all $\mathbb S^N$.

Let $T$ denote the linear integral operator defined by the right hand side of \eqref{eh}.
Then $T$ is a self-adjoint compact operator on $L^2(\mathbb S^N)$, whose spectral decomposition can be described in terms of the spaces $H_l$, $l\geq 0$ that consists of restrictions to $\mathbb S^N$ of homogeneous harmonic polynomials in $\R^{N+1}$ of degree $l$. Then $L^2(\mathbb S^N)$ is the closure of the direct sum of $H_l$, $l\ge 0$, and the elements in $H_l$ are eigenvectors of $T$ with eigenvalue
\begin{align}
\label{ee}
e_l =
\kappa_N 2^{\alpha} (-1)^l \frac{\Gamma(1-\alpha) \Gamma(N/2-\alpha)} { \Gamma(-l+1-\alpha) \Gamma(l+N-\alpha)}
\end{align}
where $\alpha = N/2-s$ and
$$
\kappa_N
=\begin{cases}
2 \pi^{1/2} & \text{if } N=1
\\
2^{2(N-1)}\pi^{(N-1)/2}\frac{\Gamma((N-1)/2)\Gamma(N/2}{(N-2)!} & \text{if } N \geq 2 .
\end{cases}
$$
Moreover  these are the only eigenvalues of $T$,
see e.g.\@ \cite{frank-lieb}.
Formula \eqref{ee} and standard identities of the Gamma function give
$$
\frac{e_{l+1}}{e_l} = \frac{l+\alpha}{l+ N-\alpha}<1 .
$$
We note that if $\phi $ is any of the functions in \eqref{kernel}, then the corresponding function $h$ belongs to $H_1$ and is a nontrivial solution of \eqref{eh}. Thus the number $a$ appearing in \eqref{eh} must be equal to $ e_1$ (this could also be deduced from an explicit formula for $a$). Therefore,  if $h$ is a solution of \eqref{eh} then  $h \in H_1$. Since the dimension of $H_1$ is $N+1$, $H_1$ is generated by functions obtained by the transformation \eqref{def h} applied to the functions in \eqref{kernel}.
The proof is concluded. \qed

\medskip

{\bf Acknowledgements.}
This work has been supported by grants Fondecyt 1090167, 1110181, CAPDE-Anillo ACT-125 and Fondo Basal CMM.


\begin{thebibliography}{00}

\bibitem{aubin}
T. Aubin, {\em Problemes isoperimétriques et espaces de Sobolev,} J. Differ. Geometry 11 (1976), 573–-598.


\bibitem{cgs}
L. Caffarelli, B. Gidas, J. Spruck, {\em Asymptotic symmetry and local behaviour of semilinear elliptic equations with critical Sobolev growth,} Commun. Pure Appl. Math. 42 (1989), 271–-297.

\bibitem{caffarelli-silvestre}
L.~Caffarelli, L.~Silvestre, {\em An extension problem related to the fractional {L}aplacian.} Comm. Partial Differential Equations 32 (2007), no. 7-9, 1245--1260.

\bibitem{cl}
E.A. Carlen, M. Loss, {\em Extremals of functionals with competing symmetries}. J. Funct. Anal. 88 (1990), no. 2, 437--456. 

\bibitem{chen-etal}
W.~Chen, C.~Li, B.~Ou, {\em Classification of solutions for an integral equation.} Comm. Pure Appl. Math. 59 (2006), no. 3, 330--343.

\bibitem{frank-lieb}
R.L. Frank, E.H. Lieb, {\em Inversion positivity and the sharp Hardy-Littlewood-Sobolev inequality.} Calc. Var. Partial Differential Equations 39 (2010), no. 1-2, 85--99. 

\bibitem{fl2}
R.L. Frank, E.H. Lieb, {\em A new, rearrangement-free proof of the sharp Hardy-Littlewood-Sobolev inequality.} In: Spectral Theory, Function Spaces and Inequalities, B. M. Brown et al. (eds.), 55 - 67, Oper. Theory Adv. Appl. 219, Birkhauser, Basel, 2012.

\bibitem{hl1}
G.H. Hardy, J.E. Littlewood, {\em Some properties of fractional integrals. I.} Math. Z. 27 (1928), no. 1, 565--606.

\bibitem{hl2}
G.H. Hardy, J.E. Littlewood, {\em On certain inequalities connected with the calculus of variations.} J. London Math. Soc. 5 (1930), 34--39.

\bibitem{talenti}
G. Talenti, {\em Best constants in Sobolev inequality,} Annali di Matematica pura e aplicata 10 (1976), 353–372
\bibitem{yanyan}
Y.Y.~Li, {\em Remark on some conformally invariant integral equations: the method of moving spheres.} J. Eur. Math. Soc. (JEMS) 6 (2004), no. 2, 153--180.

\bibitem{lz} 
 Y.Y. Li, M. Zhu, {\em Uniqueness theorems through the method of moving spheres.} Duke Math. J. 80 (1995), no. 2, 383–-417.

\bibitem{lieb}
 E. Lieb,  {\em Sharp constants in the Hardy-Littlewood-Sobolev and related inequalities.} Ann. of Math. (2) 118 (1983), no. 2, 349–-374.

\bibitem{rey}
 O. Rey {\em The role of the Green's function in a nonlinear elliptic equation involving the critical Sobolev exponent.} J. Funct. Anal. 89 (1990), no. 1, 1–-52.

\bibitem{sobolev}
S.L. Sobolev, {\em On a theorem of functional analysis} Mat. Sb. (N.S.) 4 (1938), 471--479; English transl. in Amer. Math. Soc. Transl. Ser. 2 34 (1963), 39--68.
\end{thebibliography}
\end{document}